\documentclass[11pt]{amsart}
\usepackage{geometry}                % See geometry.pdf to learn the layout options. There are lots.
\usepackage{graphicx}
\usepackage{amssymb}
\usepackage{epstopdf}
\newcommand{\R}{\mathbb{R}}
\newcommand{\T}{\mathbb{T}}
\newcommand{\eps}{\varepsilon}

\title{Normal hyperbolicity for non-autonomous oscillators and oscillator networks}
\author{R.S.MacKay}
\address{Mathematics Institute and Centre for Complexity Science,
University of Warwick, Coventry CV4 7AL, UK}
\email{R.S.MacKay@warwick.ac.uk}
\date{\today}                                           % Activate to display a given date or no date

\begin{document}
\maketitle

\section{Introduction}
Aneta Stefanovska expressed a vision ``to build a self-consistent theory of non-autonom\-ous oscillators" (June 2014).  In this direction she introduced the class of ``chronotaxic" systems \cite{SCS}, defined as ``oscillatory systems with time-varying, but stable, amplitudes and frequencies''.

This chapter presents a view of a non-autonomous oscillator as a mapping from input functions of time to a circle of possible solutions (state functions of time).  It indicates how this view encompasses chronotaxic systems and enables one, at least conceptually, to understand the extent of synchronisation in networks of oscillators, whether autonomous or not.  For the latter a hierarchical aggregation scheme is introduced.  

The approach is based on the theory of normal hyperbolicity \cite{F,HPS}.  This theory is the mathematical expression of Haken's slaving principle \cite{Ha}, the idea that some variables for a dynamical system might contract relatively rapidly onto some invariant submanifold in the state space, and then it suffices to study the dynamics on the submanifold.  Two key results of normal hyperbolicity theory are:~(i) conditions guaranteeing existence of such a submanifold, and (ii) smooth persistence of normally hyperbolic (NH) submanifolds as parameters are varied smoothly.  It was developed before Haken's slaving principle and deserves to be better known in the physics community.  It is a generalisation of centre manifold theory, which is the main mathematical tool Haken used, but has much wider scope.  An obstacle is that it demands considerable technical expertise in mathematical analysis.  Yet the obstacles are genuine:~it turns out that NH submanifolds are differentiable some number $r$ times, depending on the ratio between normal and tangential contraction rates, but typically not more than $r$ times.  This is important to recognise, as there is a tendency in physics to consider such functions as pathologies (though physicists do understand that there can be fractal functions).

It is a project on which I have been working for many years, notably with PhD student Stephen Gin (2006--13).  It was prompted initially by Mohammad Ghaffari Saadat in 2003, who had formulated a limit-cycle model for a bipedal robot walking down a slope \cite{TGN} and asked me how much non-uniformity of slope it could cope with.  I proposed to tackle this problem by fitting it into the framework of the non-autonomous version of the theory of NH submanifolds, where the result of a not too large forcing function on an oscillator is a circle of possible trajectories.   Gin and I attempted to develop good versions of the proofs of normal hyperbolicity results to produce realistic conditions guaranteeing the outcome \cite{G}.  Our approach is still incomplete, but I present here the key ideas.

%After giving a colloquium on the topic in Groningen in 2013, I discovered that my approach fits to some extent in a framework that had been proposed by Willems \cite{W}, but it seems he had not taken into account the potential use of normal hyperbolicity theory.

In the world of conservative dynamics, an oscillator is considered to be a Hamiltonian system with an elliptic equilibrium point; this view has fundamental importance but is not the appropriate one for present purposes.

Outside the world of conservative dynamics, an oscillator is usually considered to be an autonomous dynamical system with an attracting periodic orbit.  The concept has been extended to cater for chaotic oscillators, but I will postpone treating that extension until near the end of this chapter.

This concept of oscillator as a system with an attracting limit-cycle, however, fails to include the many situations where it is subject to time-dependent forcing.  Also, in a network of oscillators, each is subject to input from others, in general time-dependent, so even if the network is autonomous it is useful to consider time-dependent forcing on each of its oscillators.

So I propose a view of an oscillator as a mapping from input functions $f$ of time to a circle's worth of solutions for its state $x$ as a function of time.  Each input function $f$ (possibly with more than one component) causes a response $x_\theta$ (a function of time) with a phase $\theta \in S^1$ labelling the different possible responses.
This view is justified by the theory of normal hyperbolicity, at least for not too strong forcing.  It is also my interpretation of chronotaxic systems.

The idea is to consider a non-autonomous system $\dot{x} = v(x,t)$ on a state space $X$ as an autonomous system in the extended state space $X \times \R$, with the real line $\R$ representing time.  The dynamics has the form
\begin{eqnarray}
\dot{x} &=& v(x,s) \label{eq:sys}\\
\dot{s} &=& 1. \nonumber
\end{eqnarray}
First suppose the vector field $v = v_0$ is independent of $s$ and $\dot{x}=v_0(x)$ has a limit cycle $\gamma$ (in the strong sense of a periodic orbit with no Floquet multipliers\footnote{The Floquet multipliers of a periodic orbit are the eigenvalues of the derivative of the return map to a transverse section.} on the unit circle).  The most relevant case for applications might be the attracting case (all Floquet multipliers inside the unit circle), but one can allow the more general situation.
Then in $X\times \R$, the extended system (\ref{eq:sys}) has an extended verson of $\gamma$, namely an invariant cylinder $\gamma \times \R$.  The trajectories form helices on the cylinder, representing the same periodic solution but shifted in $s$.  

This cylinder is an example of a NH submanifold.  In general, a {\em NH submanifold} for a $C^1$ dynamical system is an invariant $C^1$ submanifold for which the linearised normal dynamics decomposes into components which contract exponentially in forward or backward time respectively, and faster than the linearised tangential dynamics.  Note that the use of the word ``normal'' might suggest perpendicular, but actually, a normal vector to a submanifold is defined to be an equivalence class of vectors at a point modulo vectors tangent to the submanifold at that point.
In the above case, the linearised tangential dynamics neither contracts nor expands on average, because the phase difference between any pair of the helices remains constant.  The linearised normal dynamics decomposes into exponentially contracting components in forward and backward time, corresponding to the Floquet multipliers inside and outside the unit circle, respectively.

Now allow $v$ to depend weakly on $s$.  The key result for NH submanifolds is that they persist under $C^1$-small  perturbation.  Thus the perturbed system has a $C^1$-nearby invariant cylinder, no longer in general of product form but diffeomorphic to $S^1\times\R$.  Furthermore, the vector field on it is close to that on the unperturbed cylinder, and the normal dynamics is close to that for the unperturbed case.  The solutions on the perturbed cylinder are not in general just a family of periodic solutions differing by phase.   In particular, there may be solutions on the cylinder to which all nearby ones converge in forward time.  There may also be solutions to which all nearby ones converge in backward time.  Or neither may happen.  In any case, there is a circle's worth of solutions on the cylinder, which one could label by the intersections of the cylinder with $s=0$ for example.  

In particular, if $v(x,t) = v_0(x) + f(t)$ then the forcing function $f$ produces a circle's worth of state functions $x$ of time on the cylinder.  In general a forcing function $f$ should be allowed to depend on the state $x$ too, so $v=v_0(x)+f(x,t)$, and by normal hyperbolicity theory, the same conclusion holds.

As an illustration, consider a model of a quasiperiodically forced limit-cycle oscillator from \cite{CS}:
\begin{eqnarray}
\dot{x}&=&-qx-\omega y \\
\dot{y}&=&\omega x - q y + \gamma f(t), \nonumber
\end{eqnarray}
with $q= \alpha(\sqrt{x^2+y^2}-a)$, $f(t) = \sin 2\pi t + \sin 4t$, $\alpha, a >0, \gamma \ge 0$ (more natural would be $q=\alpha(x^2+y^2-a^2)$ because it makes the dynamics smooth at the origin, but the interest is in the behaviour for $r=\sqrt{x^2+y^2}$ near $a$).
In polar coordinates $(r,\theta)$ and extended state-space, this is
\begin{eqnarray}
\dot{r}&=&-\alpha r(r-a) + \gamma f(s) \sin \theta \\
\dot{\theta}&=& \omega-\frac{\gamma}{r}f(s)\cos\theta \nonumber \\
\dot{s}&=&1. \nonumber
\end{eqnarray}
For $\gamma=0$ there is an invariant cylinder $r=a$.  It attracts exponentially with exponent $-\alpha a$ and the motion on the cylinder is $\dot{\theta}=\omega$, $\dot{s}=1$, which has Lyapunov exponents 0.  So the cylinder is NH and persists to a deformed invariant cylinder for small enough $\gamma$.  A rough estimate of the range of $\gamma$ for which persistence is guaranteed is given by the range for which tangential contraction is weaker than normal contraction on the unperturbed cylinder.  The normal contraction rate (onto the unperturbed cylinder) is still $\alpha a$.  The tangential contraction (or expansion if negative) $-\frac{\partial \dot{\theta}}{\partial \theta} = -\frac{\gamma}{r}f(s)\sin\theta$.  This is smaller than $\alpha a$ for all $s, \theta$ iff $2\gamma < \alpha a^2$.  Thus one can expect the NH cylinder to persist for $\gamma$ up to something of the order of $\alpha a^2/2$.  

When $\gamma$ exceeds $\alpha a^2/2$ one can not expect the invariant cylinder to persist.  It is shown numerically in \cite{CS} that the cylinder is replaced by a (non-autonomous) chaotic attractor with one unstable Lyapunov exponent (coming from the $s,\theta$ for which the tangential dynamics is expanding).  For a class of examples where a NH submanifold (in fact two 2-tori) can be proved to break up, see \cite{BaM1}.  In this chapter, however, I will concentrate on regimes of weak enough coupling that NH submanifolds persist.

As an aside, this view of an oscillator fits in Willems' ``behavioural approach'' to systems and control \cite{W}.  His view was that the description of a dynamical system should be considered to be the restrictions on the set of possible functions of time for all variables.  Normal hyperbolicity strikes me a key tool for delivering his approach.
On the other hand, he also proposed that one should go beyond the idealisation of inputs and outputs by treating all coupling as two-way, a line that I shall not follow consistently.

In this chapter I will explain how this view of an oscillator illuminates the phenomena of phase-locking, synchronisation and chimera \cite{AS}, allows to extend the concept of coupling, and allows a hierarchical reduction treatment of synchronisation in networks of oscillators.  I will extend the results to allow excitable oscillators and chaotic oscillators.  I will outline how the theory of normal hyperbolicity underlies the results.  There is a huge literature on synchronisation, e.g.~\cite{PRK}, and much of what I will say will be familiar but the important emphasis here is on synchronisation in aperiodically forced systems, which has been treated much less.

Perhaps this direction is not what Aneta had in mind, but I believe it provides a self-consistent theory for non-autonomous oscillators and I hope that it will be useful.

\section{Phase-locking}
\label{sec:pl}
It is well-known that an oscillator may phase-lock to some features of its inputs.  Indeed, this is the principle of phase-locked loops in electronic engineering \cite{Br} and of synchronous generators and motors in AC electrical networks.  

My definition of phase-locking of an oscillator to forcing is that the NH cylinder (assumed attracting) has an attracting trajectory on it and the initial condition is in its basin of attraction.  

Any discussion of attractors for non-autonomous systems requires care because the dynamics is unbounded in the time direction of extended state-space, so there are inequivalent choices of neighbourhoods of a trajectory.  For example, for the 2D system $\dot{x}=x,\dot{s}=1$, any trajectory has a neighbourhood of attraction, despite looking unstable, e.g.~for the solution $x=0$ just take neighbourhood of the form $|x|<\eps e^{2s}$.
So I make precise here that by ``attracting trajectory'' I mean the case with zero unstable space of a uniformly hyperbolic trajectory in the non-autonomous sense.  To explain what this means would take some space, so I refer the reader to \cite{BiM} (with my PhD student Zahir Bishnani), but the important feature is to choose a notion of distance in extended state-space that is uniform in time (so that one does not allow neighbourhoods like that in the above example).  There might be a bundle of trajectories which all converge together in forward time, but in general there is only one trajectory in the bundle that has a uniform-in-time neighbourhood of attraction.  It is a pullback attractor (for this concept, see the contribution by Kloeden in this volume).  My concept of attracting trajectory is distinct, however, from that of pullback attractor, because it can also occur that a pullback attractor is not uniformly hyperbolic (it may be repelling after some time).

An alternative way to describe phase-locking is that the oscillator is synchronised to its inputs.  I use ``synchronise" in a weak sense: that to a given input function of time there is a locally unique forwards asymptotic solution (the strong sense applies to systems of identical oscillators with a symmetry of the coupling that maps any oscillator to any other, and consists in all oscillators doing the same; for an example, see \cite{YM}).  Note that a forced oscillator may have more than one such attracting trajectory; this would allow different synchronisations to the same input.

This is in contrast to non-synchronisation, where there is a circle's worth of solutions that do not converge asymptotically to a discrete subset.  The strongest version of non-synchronisation is when there is a time-dependent choice of $C^1$ coordinate $\phi$ around the cylinder, replacing an initial coordinate $\theta$, such that $\dot{\phi}=\omega(t)$, a positive function of $t$ only, and $\frac{\partial \phi}{\partial \theta}$ and its inverse are bounded.  Then with a new time $\tau$ defined by $d\tau/dt=\omega(t)$, we obtain $d\phi/d\tau = 1$.  It would be interesting to investigate the probability of this case with respect to a distribution of oscillator frequencies for given weak forcing, perhaps obtaining a sort of non-autonomous KAM result\footnote{The original KAM theory gives a set of invariant tori for near-integrable Hamiltonian systems, the measure of whose complement goes to zero as the perturbation from integrability goes to zero.}, extending the theory of reducibility of cocycles (see \cite{DS} for an early example).

The main conclusion of this section is that synchronisation of an oscillator to its inputs is dimension-reduction.  In particular, if there is no immediate feedback from the oscillator to any of its inputs, then one could delete that oscillator, replacing its outputs by some modifications of the outputs from its inputs. 

\section{Synchronisation of two oscillators}
\label{sec:2osc}
Let us start with two autonomous oscillators $x_i = v_i(x_i)$, $i=1,2$, meaning each has a limit cycle $\gamma_i$, and couple them in the standard sense of a modification to the vector field of the product system, depending on the state of each but not too strongly, so 
\begin{equation}
\dot{x}_i = v_i(x_i) + g_i(x_1,x_2),
\end{equation}
with $g_i$ $C^1$-small.
Then the product system has a NH 2-torus, being a small perturbation of $\gamma_1 \times \gamma_2$.

If the difference of the frequencies of the uncoupled limit cycles is smaller in a suitable dimensionless sense than the coupling then the NH torus has an attracting limit cycle on it, which makes one turn in the $\gamma_2$ direction for each turn in the $\gamma_1$ direction.  I say the two oscillators have gone into $1:1$ synchronisation.  Recall Huygens' clocks.  The torus may have more than one attracting limit cycle on it, in which case several synchronised solutions are possible.  It may also have unstable limit cycles on it.

Similarly, if the frequencies are close to being in (coprime) integer ratio $m:n$ then coupling might produce an attracting $m:n$ limit cycle on the NH torus, which makes $m$ revolutions in the $\gamma_1$ direction and $n$ in the $\gamma_2$ direction per period.  On the other hand, for weak coupling and smooth enough dynamics, the non-synchronised situation occurs with high probability.  More precisely, if one adds a free parameter varying the unperturbed frequency ratio, then KAM theory gives a set of parameter values of nearly full measure for which the dynamics is conjugate to a constant vector field on a 2-torus with irrational frequency ratio (e.g.~\cite{LD} for a version by my PhD student Jo\~ao Lopes Dias).  Thus synchronisation does not always result.

Now consider the non-autonomous situation, where one or both of the oscillators is subject to external forcing.  If the forcing is not too strong then the resulting system has a NH submanifold in extended state space, diffeomorphic to $\gamma_1 \times \gamma_2 \times \R$, which I call a torus-cylinder.  More generally, for any manifold $M$ I define an $M$-cylinder to be a manifold diffeomorphic to $M\times\R$.  Thus an ordinary cylinder can be called a circle-cylinder.
If the unperturbed frequencies are close to integer ratio $m:n$ then the NH submanifold might contain a NH attracting submanifold diffeomorphic to a circle cross time, being a perturbation of the product of a $m:n$ synchronised limit cycle for the autonomous system and time.  In this situation the non-autonomous pair of oscillators can be replaced by a single one.

So again, synchronisation of two oscillators is a dimension-reduction.

\section{What is coupling?}
In the previous section I used the standard dynamical systems notion for coupling as a perturbation of the product of two vector fields.  One might want, however, to allow more general forms of coupling, for example incorporating time-delays or coupling via an intermediate dynamical system.  Furthermore, suppose one achieved a dimension-reduction as in section~\ref{sec:pl} or \ref{sec:2osc} and then wants to consider how the new effective oscillator is coupled to others that originally were coupled to one or both of the pair of oscillators.  This is no longer describable as a standard perturbation of the product of vector fields.

So I generalise the notion of coupling of two non-autonomous oscillators.  As already defined, a non-autonomous oscillator is a non-autonomous system with NH cylinder on which the dynamics can be described by one phase $\theta$ with $\dot{\theta} = f(\theta,t)$.  A coupling of two non-autonomous oscillators is a non-autonomous system with a NH torus-cylinder on which the dynamics can be described by two phases $\theta = (\theta_1,\theta_2)$ with $\dot{\theta}_i=\tilde{f}_i(\theta,t)$ and $\tilde{f}_i(\theta,t)$ close to $f_i(\theta_i,t)$ for some $f_i$.  

Then the dynamics on the NH torus-cylinder may contain a NH attracting circle-cylinder, as in the more restricted case of the previous section.  If the trajectory is in its basin of attraction, I say the two oscillators synchronise.

\section{Synchronisation of $N$ oscillators}
Not too strong coupling of $N$ non-autonomous oscillators produces a NH $N$-torus-cylinder.  The dynamics on it might contain an attracting NH $d$-torus-cylinder for some $d<N$.  If $d=1$ the whole group is synchronised and can be replaced by a single effective non-autonomous oscillator.  If $d=0$ the whole group is phase-locked to its inputs and can be eliminated.

Once again, synchronisation, whether partial or complete, means dimension-reduction.

\section{Hierarchical aggregation}
In a network of oscillators, the above dimension-reductions can in principle be iterated.  First one identifies groups of oscillators which synchronise or phase-lock to their inputs.  One reduces to a new network of effective oscillators.  Then one repeats, if possible.  The end result is a decomposition into synchronised clusters.

Although I did not find out about his work until after I'd proposed this, it is a direct example of Willems' ``tearing, zooming, linking'' approach \cite{W}.

One should note that the end result is not necessarily complete synchronisation.  Indeed, it could well be a chimera \cite{AS}, meaning a system in which some of the oscillators are synchronised but others behave chaotically.  The chaotic ones force the synchronised ones and the synchronised ones force the chaotic ones, but our approach of non-autonomous oscillators caters for both of these.  There is now a huge literature on chimera.  To me the phenomenon was not a surprise because it fits in my framework, but without the framework it can admittedly be considered surprising.

\section{Normal hyperbolicity estimates}
To achieve the above dimension-reductions requires good normal hyperbolicity estimates, i.e.~results guaranteeing existence of NH submanifolds.

The easiest case, namely, 1D submanifolds, which are just uniformly hyperbolic trajectories of non-autonomous systems, was already treated in \cite{BiM} (incidentally, it was formulated with attracting trajectories in mind, but another application would be to the unstable trajectories of geophysical flows that form boundaries between trajectories of different classes, e.g.~\cite{FPET}).  So that takes care of the case of phase-locking.  

Higher-dimensional NH submanifolds, however, require more theory.  The classic references are \cite{F,HPS}.  They are not particularly well adapted to producing practical estimates.  Thus I set Stephen Gin onto developing a better way.  His PhD thesis \cite{G} gives the outcome, but it is not a complete treatment.  So here, I sketch an approach to NH estimates that I believe will be useful.  It is in the classic dynamical systems setting of a vector field on the product of state space and time, but hopefully could be extended to take care of the more general forms of coupling that I have described here.

I restrict attention to submanifolds that are torus-cylinders, but of arbitrary dimension $m+1$.  
So suppose
\begin{eqnarray}
\dot{\theta} &=& \Theta(\theta,r,t) \\
\dot{r} &=& R(\theta,r,t), \nonumber
\end{eqnarray}
for $\theta \in \T^m$, $r \in U$, a neighbourhood of $0\in \R^p$.
I suppose that the product $|R_\theta| |\Theta_r|$ is small (where subscript denotes derivatives), the $r$-dynamics is hyperbolic, and the Green function for linearised normal dynamics decays faster than any contraction that may occur in $\theta$-dynamics.

Given a Lipschitz graph $r=\rho(\theta,t)$, a candidate for an invariant submanifold, construct a new one, $T\rho$, by the following steps:
\begin{enumerate}
\item For all $(\theta_0,t_0)$, let $\theta()$ be the trajectory of $\dot{\theta}(t)=\Theta(\theta,\rho(\theta,t),t)$ from $\theta(t_0)=\theta_0$.
\item Solve $\dot{r}(t) = R(\theta(t),r(t),t)$ for the unique function $r()$ such that $r(t)$ is near $\rho(\theta(t),t)$ for all $t$.
\item Set $(T\rho)(\theta_0,t_0) = r(t_0)$.
\end{enumerate}
To achieve the second step, I assume that $L:C^1(\R,\R^{p})\to C^0(\R,\R^{p})$ defined by
$$L[x](t) = \dot{x}(t) - R_r(\theta(t),r(t),t) x(t)$$
on infinitesimal displacements $x$ in $r$ has bounded inverse.  This is equivalent to the first part of the NH condition, namely a splitting of the normal bundle into exponentially contracting backwards and forwards subspaces.

Having thus constructed the ``graph transform'' $T$, I want to prove that it is a contraction on a suitable space of graphs and hence has a unique fixed point there, which will be an invariant graph.  In the direction of achieving this,
define a {\em slope} to be a linear map $\sigma$ from displacements in $\theta$ to displacements in $r$.
For an approximation $\tilde{\sigma}$ to the expected derivative $\rho_\theta$, define
$M_{\tilde{\sigma}}: W^{1,\infty}(\R,\R^{mp}) \to W^{0,\infty}(\R,\R^{mp})$ by
$$M_{\tilde{\sigma}}[\sigma] = \dot{\sigma} - R_r\sigma+\sigma(\Theta_\theta + \Theta_r\tilde{\sigma})$$
on slope functions $\sigma$ of $t$, where $W^{s,\infty}$ are the spaces of functions with essentially bounded $s^{th}$ derivative.
Suppose that $M_{\tilde{\sigma}}$ has bounded inverse.  This is the second part of the NH condition, namely faster normal contraction than tangential contraction.

Then $T$ should be a contraction in the space of $C^0$ functions with an a priori Lipschitz constant.  So it would have a unique fixed point $\rho$.  Any fixed point is invariant and actually $C^1$ with slope $\rho_\theta$ being the fixed point of the contraction map $\sigma \mapsto M_\sigma^{-1}[R_\theta]$.  

To complete this programme requires detailed estimates.  Formulated in terms of contraction maps as here, it should be possible to obtain excellent estimates, along the lines of the uniformly hyperbolic case in \cite{BiM}.  We might do best to follow the approach of \cite{H} (cf.~\cite{E}), but replacing their exponential hypotheses by our hypotheses of invertibility of $L$ and $M$ and modifying their exponentially weighted norm to use the linearised tangential flow.  I would like to finish this one day.

\section{Extension to class 1 neurons}
So far, I have considered the simplest type of oscillator, namely limit cycles, but the treatment can be extended to class I neurons (or excitable oscillators).  These are dynamical systems with an attracting invariant cylinder in the autonomous case and dynamics on it in simplest form given by
\begin{eqnarray}
\dot{\theta} &=& \mu + 1-\cos\theta \\
\dot{\mu}&=&0 . \nonumber
\end{eqnarray}
Since $\mu$ is constant, one could think of it as an external parameter, but I wish to consider it as a state variable because coupling from another neuron can make $\mu$ change in time.
It is best to think of $\mu$ as bounded, so the attracting cylinder can be considered an invariant annulus.

They arise in modelling of ``excitable neurons'' whose frequency goes to zero as a parameter ($\mu$) is varied and then settle at a $\mu$-dependent resting state, or in reverse go from a resting state to large amplitude periodic spiking.
An example is the Morris-Lecar model \cite{ML}, but it was \cite{EK} who identified the phenomenon as the unfolding of a saddle-node on a cycle (I proposed this independently to physiologist H.Barlow in the same year and then in 1991 proposed to C.Koch the extension to allow crossover at a ``saddle-node loop'' \cite{S} to the unfolding of a homoclinic orbit to a saddle).
Thus the non-autonomous version has an attracting NH annulus-cylinder.

I had an undergraduate student study networks of such neurons in 1989/90, with the state $\mu$ of each neuron driven by the spiking of some others (with time-delay kernels), which produced periodic bursting \cite{M2}.

Two class I neurons coupled not too strongly have a NH attracting annulus$\times$annulus-cylinder.  Generic bifurcation diagrams in the autonomous case were given in \cite{BaM}.  The dynamics on it has attracting submanifolds of various types.  The non-autonomous case has non-autonomous versions of them.

The theory of this paper applies just as well to class I neurons as to ordinary oscillators, with the addition of the $\mu$-direction for each class I neuron.

\section{Extension to chaotic oscillators}
The approach can also be extended to chaotic oscillators if they have an attracting NH submanifold containing the attractor.  For example, think of a R\"ossler attractor \cite{R}, which is contained in a solid torus in $\R^3$.  Then the non-autonomous system has a solid-torus-cylinder.  A R\"ossler attractor can be phase-locked to forcing, meaning that the dynamics is attracted onto a disk-cylinder (a solid torus is the product of a disk and a circle).  This should be quite easy because the R\"ossler attractor was observed to be nearly phase-coherent.  I interpret that as meaning that there is a cross-section with nearly constant return time (equivalently, for a given cross-section $\Sigma$ there is a constant $c>0$ and a function $b:\Sigma \to \R$ such that the return time $\tau(x) = c + b(f(x))-b(x)$, where $f:\Sigma \to \Sigma$ is the return map).

Synchronisation of chaotic attractors with NH cylinders of dimensions $N_1+1, N_2+1$ means there is a NH cylinder for the coupled system with dimension less than $N_1+N_2+1$.

Even better, the theory of NH submanifolds extends to NH laminations \cite{HPS}.  A lamination is a topological space in which each point has a neighbourhood homeomorphic to the product of a Euclidean space with a general topological space.  It decomposes into leaves, which are locally submanifolds but in general only injectively immersed, so a leaf may accumulate onto itself.  The theory of NH laminations requires a $C^1$-structure in addition, but is basically the same as for NH submanifolds.  In particular, a NH lamination persists under $C^1$-small perturbation.  

This means one can treat some chaotic attractors in greater detail.  In particular, imagine we start with a non-trivial uniformly hyperbolic attractor of an autonomous system, for example a suspension of a Plykin attractor \cite{P}.  This is perhaps less familiar than R\"ossler's attractor but deserves to be better known, as the simplest uniformly hyperbolic attractor after equilibria and periodic orbits.  The Plykin attractor was constructed for a discrete-time system, but the map is isotopic to the identity so one can realise it as the first return map of an associated continuous-time system.  My PhD student Tim Hunt showed an explicit way to realise it in a system of three ODEs, extended by another PhD student Linling Ru, and less cumbersome ways have been proposed (though not yet with rigorous justification) \cite{K}.  It is a NH lamination, whose leaves are its unstable manifolds (of dimension two:~one expanding dimension and one time dimension) and they form a Cantor set transversally.  Under time-dependent forcing, it persists to a Cantor set of 3D leaves whose tangent space is spanned by one expanding dimension and two near neutral dimensions.  The persistence is highly robust, requiring only that any tangential contraction be slower than any transverse contraction.  Then one can ask what happens on the leaves.  The dynamics might collapse onto a 2D subleaf with the same expanding dimension one neutral dimension.  I would say the attractor has synchronised to the forcing.

Similarly, one could couple a suspended Plykin attractor to a limit-cycle oscillator.  It produces an attractor with a Cantor set of 3D leaves (the product of the 2D leaves of the chaotic attractor with the limit cycle).  The dynamics of each leaf might collapse onto 2D subleaves.  I would say the Plykin attractor and limit cycle synchronise together.

More generally, one could couple a continuous-time autonomous uniformly hyperbolic attractor with $M$ unstable dimensions to $N$ limit cycle oscillators and obtain an attractor with a mixture of chaos and nearly quasiperiodic behaviour.  It would have $M$ unstable dimensions, $N$ nearly quasiperiodic dimensions, and the flow dimension, with the remaining dimensions contracting onto the leaves.  By the theory of NH laminations, such attractors persist for small smooth perturbations, though the dynamics in the quasiperiodic dimensions cannot be expected to remain quasiperiodic.  Nonetheless, it will have small Lyapunov exponents for those dimensions and perhaps there is a non-autonomous KAM theory that would even give truly quasiperiodic motion for a set of nearly full measure of parameters.  I propose this as an explanation of the scenario reported recently by \cite{YK}. 

As a final note, one might ask about physical realisation of attractors like R\"ossler's.  I designed an electronic oscillator back in 1981, principally to demonstrate period-doubling sequences \cite{M1}, but moving the parameter further it exhibited a R\"ossler type of attractor.  Model equations for the voltages at three points have the form
\begin{eqnarray}
\dot{x}&=&ax-by \\
\dot{y}&=& cx-ez \nonumber \\
\dot{z} &=& -fy-g(z) , \nonumber
\end{eqnarray}
with $a,b,c,e,f$ positive constants of which $a$ was adjustable by a 10-turn potentiometer, and $g$ an approximately odd cubic nonlinearity produced with a pair of transistors.  Interestingly, as I increased $a$ further, the R\"ossler attractor turned into what Chua later called a double-scroll attractor \cite{MCK}.  Indeed, Chua's equations turn out to be equivalent to mine after minor changes of variable.
 
\section{Conclusion}
I have shown that the behaviour of networks of oscillators, autonomous or not, can be aided by identifying normally hyperbolic submanifolds.  This allows a deeper understanding of synchronisation of oscillators to forcing and to each other, especially in the aperiodic case.  There are many studies on synchronisation in autonomous or periodically forced systems (for one example, see \cite{SST}) but relatively few on the aperiodically forced case.  The fundamental feature of synchronisation is dimension-reduction of an associated normally hyperbolic submanifold.  In a network of oscillators, even if autonomous, the inputs that an individual oscillator sees are in general aperiodic.  This motivates a hierarchical aggregation scheme for understanding the dynamics of a network of oscillators:~oscillators that synchronise to their inputs can be eliminated, groups of oscillators that synchronise together can be replaced by a single effective oscillator.  All this depends on generalising the notion of oscillator from a limit cycle of an autonomous dynamical systems to a mapping from input functions of time to solutions and generalising the notion of coupling.  Finally, I extended the treatment from limit-cycle oscillators to excitable oscillators and chaotic oscillators.

\end{document}